\documentstyle[12pt]{article}
\newcommand{\const}{\mathop{\rm const}\limits}

\newcommand{\vraisup}{\mathop{\rm vraisup}\limits}

\begin{document}

\begin{center}

{\bf EXACT EXPONENTIAL TAIL ESTIMATIONS }\\

\vspace{4mm}

{\bf IN THE LAW OF ITERATED LOGARITHM  }\\

\vspace{4mm}

{\bf FOR BOCHNER'S MIXED LEBESGUE SPACES }\par

\vspace{4mm}

 $ {\bf E.Ostrovsky^a, \ \ L.Sirota^b } $ \\

\vspace{4mm}

$ ^a $ Corresponding Author. Department of Mathematics and computer science, Bar-Ilan University, 84105, Ramat Gan, Israel.\\
\end{center}
E - mail: \ galo@list.ru \  eugostrovsky@list.ru\\
\begin{center}
$ ^b $  Department of Mathematics and computer science. Bar-Ilan University,
84105, Ramat Gan, Israel.\\

E - mail: \ sirota3@bezeqint.net \\

\vspace{4mm}
                    {\sc Abstract.}\\

 \end{center}

 \vspace{4mm}

   We  obtain the quite exact exponential bounds for tails of distributions of sums of Banach space valued random variables
uniformly over the number of summands  under natural  for the Law of Iterated Logarithm  (LIL) norming. \par

  We  study especially  the case of the so-called  mixed (anisotropic)  Lebesgue-Riesz spaces, on the other words,
 Bochner's spaces, for instance, continuous-Lebesgue spaces,
 which appear for example in the investigation of non-linear Partial Differential Equations of evolutionary type.\par

  We give also some  examples in order to show the exactness of our estimates.\par

  \vspace{4mm}


\vspace{4mm}

 \par

\vspace{4mm}

\section{Notations. Statement of problem.}

\vspace{3mm}

{\bf 1.}  Let  $  (B, ||\cdot||B )  $  be separable non - trivial: $ \dim(B) \ge 1 $
Banach space and  $ \{  \xi_j   \}, \xi = \xi_1, \ j = 1,2,\ldots  $
be a sequence of centered in the weak sense: $ {\bf E} (\xi_i,b) = 0 \ \forall b \in B^* $  of independent identical distributed
(i. i.d.)  random variables (r.v.) (or equally random vectors, with at the same abbreviation r.v.) defined on some non-trivial
probability  space $  (\Omega = \{\omega\}, F, {\bf P})   $ with values  in the space  $ B. $  Denote

 $$
 S(n) =  \sum_{j=1}^n \xi_j, \ n = 1,2,\ldots. \eqno(1.1)
 $$

We suppose also that the r.v. $ \xi  $ has a weak second moment:

 $$
 \forall b \in B^* \ \Rightarrow  (Rb,b) := {\bf E } (\xi,b)^2 < \infty.
 $$

 Let also  $ v = v(n) $ be  arbitrary non-random  positive strictly increasing numerical sequence   such that

 $$
 v(1) = 1, \hspace{5mm} \lim_{n \to \infty} v(n) = \infty, \eqno(1.2)
 $$
for instance,

$$
v_r(n) :=  \left[\log( \log ( n + e^e -1) \ ) \right]^r,  \ r = \const \ge 1/2. \eqno(1.3)
$$

 We set

 $$
 \tau(n) = \frac{S(n)}{\sqrt{n} \ v(n)}, \hspace{5mm} \tau_r(n) = \frac{S(n)}{\sqrt{n} \ v_r(n)}, \eqno(1.4)
 $$

$$
Q(u) = Q^{(\xi)} (u) \stackrel{def}{=} {\bf P} ( \sup_n ||\tau(n)|| > u), \ \eqno(1.5)
$$

$$
 Q_r(u) = Q_r^{(\xi)} (u)  \stackrel{def}{=} {\bf P} ( \sup_n ||\tau_r(n)|| > u), \
u = \const \ge e. \eqno(1.6)
$$

\vspace{4mm}

{\bf  Our purpose in this article is obtaining some sufficient conditions for exponential decreasing as $ u \to \infty $
of the probabilities $ Q(u), \ Q_r(u),  $ for example, of a view }

 $$
 Q(u) \le \exp \left(-C u^{\beta_1} \ \log^{\beta_2}(u) \right), \ u > e, \ C,\beta_{1,2} = \const, \ C,\beta_1 > 0. \eqno(1.7)
 $$

\vspace{4mm}

{\bf 2.} The problem of describing of necessary (sufficient) conditions  for the infinite - dimensional LIL and closely
related CLT in Banach space $  B $
has a long history; see, for instance, the  monographs \cite{Araujo1} - \cite{Ostrovsky1} and articles
\cite{Garling1} - \cite{Zinn1} (CLT); \cite{Acosta1} - \cite{Ostrovsky403} (LIL);
see also reference therein.\par
 The case when $ B  = C(T), $ where $  T  $ is metrizable compact set, is considered in fact in  \cite{Ostrovsky403},
 see also \cite{Ostrovsky1}, chapter 2, section 2.7. \par
  The applications of considered theorem in statistics and method Monte-Carlo  see, e.g. in
\cite{Frolov1} -  \cite{Talay1}.\par

\vspace{3mm}

{\bf 3.} The {\it lower bound} for the probability $  Q(u) $ is very simple: as long as $ v(1) = 1, $

$$
Q(u) \ge {\bf P} (\tau(1) > u) = {\bf P} (||\xi|| > u). \eqno(1.8)
$$
 Therefore, if we want to establish the inequality (1.7), we must adopt the condition

$$
{\bf P} (||\xi|| > u) \le \exp \left(-C u^{\beta_1} \ \log^{\beta_2}(u) \right), \ u > e, \ C,\beta_{1,2} = \const, \ C,\beta_1 > 0. \eqno(1.9)
$$

\vspace{3mm}

{\bf 4.} We recall here the definition and some important properties of
 the so-called anisotropic Lebesgue (Lebesgue-Riesz) spaces, which appeared in the famous article of
 Benedek A. and Panzone  R. \cite{Benedek1}.  More detail information about this
spaces with described applications see in the books  of  Besov O.V., Il’in V.P., Nikol’skii S.M.
\cite{Besov1}, chapter 1,2; Leoni G. \cite{Leoni1}, chapter 11; \cite{Lieb1}, chapter 6. \par

\vspace{3mm}

  Let $ (X_k,A_k,\mu_k), \ k = 1,2,\ldots,l $ be measurable spaces with sigma-finite {\it separable}
non - trivial measures $ \mu_k. $ The separability  denotes that  the metric space
$ A_k  $  relative the distance

$$
\rho_k(D_1, D_2) = \mu_k(D_1 \Delta D_2) = \mu_k(D_1 \setminus D_2) + \mu_k(D_2 \setminus D_1) \eqno(1.10)
$$
is separable.\par
Let also $ p = (p_1, p_2, . . . , p_l) $ be $ l- $ dimensional vector such that
$ 1 \le p_j < \infty.$ \par

 Recall that the anisotropic  (mixed) Lebesgue - Riesz space $ L_{ \vec{p}} $ consists on all the  total measurable
real valued function  $ f = f(x_1,x_2,\ldots, x_l) = f( \vec{x} ) $

$$
f:  \otimes_{k=1}^l X_k \to R
$$

with finite norm $ |f|_{ \vec{p} } \stackrel{def}{=} $

$$
\left( \int_{X_l} \mu_l(dx_l) \left( \int_{X_{l-1}} \mu_{l-1}(dx_{l-1}) \ldots \left( \int_{X_1}
 |f(\vec{x})|^{p_1} \mu(dx_1) \right)^{p_2/p_1 }  \ \right)^{p_3/p_2} \ldots   \right)^{1/p_l}. \eqno(1.11)
$$

 In particular, for the r.v. $ \xi  $

$$
  |\xi|_p =  \left[ {\bf E} |\xi|^p \right]^{1/p}, \ p \ge 1.
$$

 Note that in general case $ |f|_{p_1,p_2} \ne |f|_{p_2,p_1}, $
but $ |f|_{p,p} = |f|_p. $ \par

 Observe also that if $ f(x_1, x_2) = g_1(x_1) \cdot g_2(x_2) $ (condition of factorization), then
$ |f|_{p_1,p_2} = |g_1|_{p_1} \cdot |g_2|_{p_2}, $ (formula of factorization). \par

 Note that under conditions of separability (1.6) on the measures $ \{  \mu_k \} $   this spaces are also  separable Banach spaces. \par

  These spaces arise in the Theory of Approximation, Functional Analysis, theory of Partial Differential Equations,
theory of Random Processes etc. \par

\vspace{3mm}

{\bf 5.} Let for example $  l = 2; $ we agree to rewrite for clarity the expression for $ |f|_{p_1, p_2}  $ as follows:

$$
|f|_{p_1, p_2} := |  f|_{p_1, X_1; p_2, X_2}.
$$
 Analogously,

$$
|f|_{p_1, p_2,p_3} = |  f|_{p_1, X_1; p_2, X_2; p_3, X_3}. \eqno(1.12)
$$

  Let us give an example. Let $ \eta = \eta(x, \omega )  $ be bi-measurable  random field, $ (X = \{x\}, A, \mu) $ be measurable space,
  $  p = \const \in [1,\infty). $ As long as the expectation $ {\bf E } $ is also an integral, we deduce

  $$
{\bf E} |\eta(\cdot, \cdot)|^p_{p,X} = {\bf E} \int_X |\eta(x, \cdot)|^p \ \mu(dx)    =
  $$

$$
 \int_X  {\bf E}  |\eta(x, \cdot)|^p \ \mu(dx)  = \int_X \mu(dx) \left[ \int_{\Omega} |\eta(x, \omega)|^p \ {\bf P}(d \omega)  \right];
$$

$$
| \ |\eta|_{p,X} \ |_{p, \Omega} = \left\{ \int_X \mu(dx) \left[ \int_{\Omega} |\eta(x, \omega)|^p \ {\bf P}(d \omega)  \right] \right\}^{1/p} =
|  \ \eta(\cdot, \cdot)  \  |_{1,X; p,\Omega}. \eqno(1.13)
$$

\vspace{3mm}

{\bf 6.}  Constants of Rosenthal-Dharmadhikari-Jogdeo- ...\par

 Let  $  p = \const \ge 1, \hspace{4mm}  \{ \zeta_k \} $ be a sequence of numerical centered, i. i.d. r.v.  with finite $ p^{th} $ moment
 $ | \zeta|_p < \infty. $  The following constant,  more precisely, function on $ p, $ is called
 constants of Rosenthal-Dharmadhikari-Jogdeo-Johnson-Schechtman-Zinn-Latala-Ibragimov-Pinelis-Sharachmedov-Talagrand-Utev...:

$$
K_R(p) \stackrel{def}{=} \sup_{n \ge 1} \sup_{ \{\zeta_k\} } \left[ \frac{|n^{-1/2} \sum_{k=1}^n \zeta_k|_p}{|\zeta_1|_p} \right].  \eqno(1.14)
$$
 We will  use  the following ultimate up to an error value $ 0.5\cdot 10^{-5} $  for the constant $ C_R $
 estimate  for $ K_R(p), $ see \cite{Ostrovsky502} and reference therein:

 $$
 K_R(p) \le \frac{C_R \ p}{ e \cdot \log p}, \hspace{5mm}  C_R = \const := 1.77638.  \eqno(1.15)
 $$
 Note that for the symmetrical distributed r.v. $ \zeta_k $ the constant $  C_R $ may be reduced  up to a value $ 1.53572.$\par

\vspace{3mm}

{\bf 7. }   Suppose the sequence $ \{\xi_i \} $ satisfies the LIL in the classical statement under our norming function:

$$
\overline{\lim}_{n \to \infty} \frac{||S(n)||}{\sqrt{n} \ v(n)} =: \zeta = \zeta(\omega),
$$
where $ {\bf P} (0 < \zeta < \infty) = 1 $  and wherein the r.v. $ \zeta $ may be non-random constant. Then

$$
\sup_n \frac{||S(n)||}{\sqrt{n} \ v(n)}  < \infty
$$
almost everywhere and following

$$
\lim_{u \to \infty} Q(u) = 0.
$$

\vspace{3mm}

{\bf 8.}  The natural norming sequence in Banach space valued LIL may be essentially greatest as the classical sequence
$ v_{1/2}(n) = \sqrt{n} \ (\log (\log( n + e^e - 1))^{1/2}. $ Namely, in the article \cite{Einmal2} there is (for any number
value $ r > 1/2) $ an example of a separable Banach space $  (B, \ || \cdot ||) $ and a sequence of i.i.d. centered r.v.
$ \{\xi_i \} $ with finite weak second moment  such that

$$
0 < \overline{\lim}_{n \to \infty} \left[ \frac{||S(n)||}{ \sqrt{n} \ (\log (\log( n + e^e - 1))^{r} )}\right] < \infty \eqno(1.16)
$$
almost surely.\par

\vspace{4mm}

\section{ Main result: exponential estimates for  LIL in ordinary Lebesgue - Riesz spaces.  }

\vspace{3mm}

  We study in this section the exponential tail estimates  for the normed sums of centered, i.i.d. r.v. in the separable Banach space
 $ L_p, \ 2 \le p < \infty, $  having the weak second moment.\par
   In detail,  $ \xi(x) =  \xi(x,\omega) $
 be bi - measurable centered  random field, $ (X = \{x\}, A, \mu) $ be  measurable space with separable sigma - finite measure $ \mu, $
 and $ \{ \xi_k(x,\omega)  \} = \{ \xi_k(x)   \} $ be independent copies of $ \xi = \xi(x) := \xi(x, \omega). $

So, in this section

$$
||\xi|| = ||\xi(\cdot, \cdot) || = |\xi|_p = \sqrt[p]{\int_X |\xi(x)|^p \mu(dx) } = |\xi|_{p,X}. \eqno(2.0)
$$

 Some additional notations. $  m = \const \ge 1, \ L = \const = p m \ge p, $

$$
\Sigma_n(x) := n^{-1/2} \sum_{k=1}^n \xi_k(x) = n^{-1/2} S(n; x), \eqno(2.1)
$$

$$
g(L) = g_p(L) = 2  \ K_R(L) \cdot \sqrt[p]{\left[ \int_X {\bf E} |\xi(x)|^L \ \mu(dx) \right]^{p/L}}, \eqno(2.2)
$$

$$
L_0 = L_0(p) = \sup \{L: \ g_p(L) < \infty \}, \eqno(2.3)
$$

$$
 h(z) = h_p(z): = \inf_{ L \in (p, L_0)} \left\{  \frac{g^L_p(L)}{z^L}  \right\}, \ z > 1. \eqno(2.4)
$$

We  clarify the role of a function $ h_p = h_p(z): $ if for the real valued r.v. $ \zeta $

$$
 \left[{\bf E} |\zeta|^L \right]^{1/L} \le g(L),
$$
then it follows from Tchebychev-Markov inequality after optimisation over $ L: $

$$
{\bf P} (|\zeta| > z) \le h(z).  \eqno(2.5)
$$

 It will be presumed hereinafter $ L_0 > p; $ may be $ L_0 = \infty; $ otherwise our propositions are trivial. \par

  Further, let $ Z_+ = (1,2,\ldots ) $ be positive integer  semi - axis; we introduce as ordinary a partition $  \Delta  $ of the set  $  Z_+: $

$$
\Delta = \{  \Delta(k) \}, \ k = 1,2,\ldots;  \  \Delta(k) = Z_+ \cap [ A(k), \ A(k+1) ), \ A(1) = 1,
$$

$$
 A(k+1) \ge A(k) + 2, \ \lim_{k \to \infty} A(k) = \infty,
$$
so that $ Z_+ = \cup_{k=1}^{\infty} \Delta(k). $\par

 Let $ w = \const > 1; $ we will say that the partition   $ \Delta = \{\Delta(k) \}  $  belongs to the class $ Y(w), $
iff

$$
\inf_k \left[ \frac{A(k+1)-1}{A(k)} \right] \ge w^2. \eqno(2.6)
$$
 For instance, $ A(k) = d^k- d +1 $  for some fixed $ \ d = 2,3,\ldots.  $

 Let at last the partition $ \Delta $ be from the set  $ Y(w), \ w > 1 $ be a given. {\it Define an important function}

 $$
 G(u) = G_{\xi}(\Delta,p,v; u) \stackrel{def}{=} \sum_{k=1}^{\infty} h_p(u \ v(A(k)) /w), \eqno(2.7)
 $$
if of course  $ G(u) $ there exists and $ \lim_{u \to \infty} G(u) = 0. $ \par

\vspace{3mm}

{\bf Theorem 2.1.  Main result.}\\

$$
Q^{(\xi)}(u) \le  G_{\xi}(\Delta,p,v; u), \ u \ge e. \eqno(2.8)
$$

{\it As a slight consequence:}

$$
Q^{(\xi)}(u) \le \inf_{w > 1} \ \inf_{\Delta \in Y(w) } \  G_{\xi}(\Delta,p,v; u), \ u \ge e. \eqno(2.9)
$$

\vspace{3mm}

{\bf Proof } used standard arguments for investigation of LIL in the (separable) Banach space.\par

\vspace{4mm}

{\bf Step 1. Lemma 2.1.} {\it  Let $ m = \const \ge 1 $ be number, not necessary to be integer,
 for which }

$$
 {\bf E} | \xi(\cdot)|_{pm,X}^{pm} < \infty. \eqno(2.10)
$$

{\it or equally }

$$
 |\xi(\cdot) \ |_{p,X; mp, \Omega}  < \infty.  \eqno(2.10a)
$$

 {\it or equally }

$$
{\bf E }  \int_X |\xi(x)|^{pm} \ \mu(dx) = \int_X {\bf E } |\xi(x)|^{pm} \ \mu(dx) < \infty.  \eqno(2.10b)
$$

{\it Then }

$$
\sup_n {\bf E} | \xi(\cdot)|_{p,X}^{pm} \le   \left[ \int_X \left[ {\bf E} | \xi(x)|^{pm}  \right]^{1/m} \ \mu(dx) \right]^m, \eqno(2.11)
$$

 {\it or equally}

 $$
 | \xi|_{p, X; mp,\Omega} \le \ |\xi|_{pm, \Omega; p,X}. \eqno(2.11a).
 $$

{\it  Recall that $ p m = L. $ }

\vspace{3mm}

{\bf Proof of Lemma 2.1.}\\

\vspace{3mm}

 Denote $  \eta(x) = |\xi(x)|^p.  $ Note that the  space $ L_m(\Omega)  $ is a Banach space; and we can apply the
generalized triangle (Minkowsky) inequality:

$$
\left[ {\bf E}  \left(  \int_X |\xi(x)|^p \ \mu(dx)  \right)^m \right]^{1/m}  =
$$

$$
\left[ {\bf E}  \left(  \int_X \eta(x) \ \mu(dx)  \right)^m \right]^{1/m}  = \left| \ \int_X \eta(x) \ \mu(dx)   \  \right|_{m, \Omega} \le
$$

$$
 \ \int_X  |\eta(x)|_{m, \Omega} \ \mu(dx) = \int_X \sqrt[m] {{\bf E} \eta^m(x) } \ \mu(dx)  =
\int_X \sqrt[m] {{\bf E} |\xi|^{pm}  (x) } \ \mu(dx), \eqno(2.12)
$$
which is equivalent to the assertions (2.11) - (2.11a). \par
 We  used theorem of Fubini-Tonelli.\par

\vspace{4mm}

{\bf Step 2.}   We apply the assertion of Lemma 2.1. to the random field
$  \Sigma_n(x)  $ instead $ \xi(x): $

$$
{\bf E} \left[ \int_X |\Sigma_n(x)|^p \ \mu(dx)  \right]^m \le
\left\{ \int_X  \left[{\bf E}|\Sigma_n(x)|^{pm} \ \right]^{1/m} \mu(dx) \right\}^m.  \eqno(2.13)
$$

 We obtain by means of Rosenthal's inequality (recall that we  consider here the case only when $  p \ge 2:) $

$$
{\bf E} |\Sigma_n(x)|^{pm} \le K_R^{p m} (p m) \ {\bf E} |\xi(x)|^{pm} =  K_R^{p m}(p m) \cdot |  \xi(x)|^{pm}_{pm, \Omega}.
$$
 It remains to substitute into (2.13):

$$
\sup_n {\bf E} | \Sigma_n(\cdot)|_{p,X}^{pm} \le K_R^{pm} (pm) \  \left[ \int_X \left[ {\bf E} | \xi(x)|^{pm}  \right]^{1/m} \ \mu(dx) \right]^m, \eqno(2.14)
$$

 or equally

 $$
 \sup_n | \Sigma_n|_{p, X; mp,\Omega} \le K_R(pm) \ |\xi|_{pm, \Omega; p,X}. \eqno(2.14a).
 $$

\vspace{3mm}

{\bf Step 3.}  Denote

$$
||\Sigma_n||^* = \max_{k=1,2,\ldots,n} ||\Sigma_k||. \eqno(2.15)
$$
 Since the sequence $ \{ \Sigma_n \}  $ under natural filtration is a  (Banach value) martingale, we can
use the Doob's inequality; the Banach space valued  martingale version is done by G.Pisier
\cite{Pisier1}, \cite{Pisier2}:

$$
| \ ||\Sigma_n||^* \ |_{L, \Omega}  \le \frac{L}{L-1} \cdot | \ ||\Sigma_n|| \ |_{L, \Omega} \le
2  \cdot | \ ||\Sigma_n|| \ |_{L, \Omega}, \eqno(2.16),
$$
since $ L \ge 2.  $\\

\vspace{3mm}

{\bf Step 4.} Now everything is ready for the final stage. Actually, let $ \Delta $ be any partition from the class $  Y(w). $
We get consequently:

$$
Q(u) = {\bf P} \left( \cup_{k=1}^{\infty} \max_{ A(k) \le n < A(k+1)} \frac{||S(n)||}{\sqrt{n} \ v(n)} > u \right) \le
$$

$$
\sum_{k=1}^{\infty} {\bf P} \left( \ \max_{ A(k) \le n < A(k+1)} \frac{||S(n)||}{\sqrt{n} \ v(n)} > u \right) \le
$$

$$
\sum_{k=1}^{\infty} {\bf P} \left( \ \max_{ A(k) \le n < A(k+1)} ||S(n)|| >  u \ \sqrt{A(k)} \ v(A(k)) \right) \le
$$

$$
\sum_{k=1}^{\infty} {\bf P} \left( \ \max_{ A(k) \le n < A(k+1)} ||\Sigma_n|| >  u \ \sqrt{A(k)/A(k+1)} \ v(A(k)) \right) \le
$$

$$
\sum_{k=1}^{\infty} {\bf P} \left( \ \max_{ A(k) \le n < A(k+1)} ||\Sigma_n|| >  u \ v(A(k))/w \right) \le
$$

$$
\sum_{k=1}^{\infty} {\bf P} \left( ||\Sigma_{A(k+1)-1 }||^* > 0.5 \ u \ v(A(k))/w \right) \le
$$

$$
 \sum_{k=1}^{\infty} h_p(u \ v(A(k)) /w) =  G_{\xi}(\Delta,p,v; u), \eqno(2.17)
$$
we used the inequality (2.5).\par
This completes the proof  of theorem 2.1.\\

\vspace{3mm}

{\bf Examples 2.1.} \\

 Suppose

 $$
 {\bf P}  \left( |\xi|_{p,X} > z \right) \le e^{-z^{\beta} }, \ z > 0, \ \beta = \const > 0;
 $$
then  of course

$$
\forall m > 0 \ \Rightarrow {\bf E} |\xi|_{p,X}^{ p m} < \infty.
$$

 The case $ \beta = \infty  $ implies the boundedness of a r.v. $ |\xi|_p:  $

$$
\vraisup_{\omega \in \Omega} |\xi|_{p,X} < \infty.
$$

 Denote

 $$
 r_0 = r_0(\beta) = \frac{\beta+1}{\beta}, \ r_0(\infty) = 1.
 $$

 We deduce after simple calculations:

$$
Q_{r_0}(u) \le \exp \left( - C(\beta,p) \ u^{\beta/(\beta+1)   }  \right),  C(\beta,p) > 0.
$$

 In particular, if $ \beta = \infty, $ then

$$
Q_{1}(u) \le \exp \left( - C \ u   \right).
$$

\vspace{3mm}

 More generally, if

$$
 {\bf P}  \left( |\xi|_{p,X} > z \right) \le e^{-z^{\beta_1} \ (\log z)^{-\beta_2} }, \ \beta_1 = \const > 0, \ \beta_2 = \const, \ z \ge e,
$$
then for some positive value $ C_3 = C_3(\beta_1, \beta_2,p) $ and $  u > e $

$$
Q_{r_0}(u) \le
\exp \left( - C_3(\beta_1, \beta_2,p) \ u^{\beta_1/(\beta_1+1)} \ ( \log u )^{( - \beta_2 -  \beta_1(\beta_1-1))/( \beta_1 + 1 ) }  \right).
$$

 We used some estimations from the monograph  \cite{Ostrovsky1}, chapter 2, section 3, p. 55 - 57; where are obtained,
in particular, the exponential tail estimates for random variables via its moment estimates.  \par

\vspace{3mm}
\section{ Exponential bounds for LIL in  Mixed Lebesgue - Riesz spaces. }
\vspace{3mm}

   Let us return to the CLT in the space  $ L_{ \vec{p}}, \ \vec{p} = \{ p_k  \}, k=1,,2,\ldots,l, \ l  \ge 2; $
 where $  1 \le p_k < \infty, $ described below. \par

 Define

$$
\overline{p} : = \max(p_1,p_2, \ldots, p_l).
$$
{\it and suppose  everywhere further}

$$
 \overline{p} \ge 2. \eqno(3.0)
$$

 Some additional notations. As in last section, $  m = \const \ge 1, \ L = \const = \overline{p} \ m \ge \overline{p}, $

$$
\Sigma_n(x) := n^{-1/2} \sum_{k=1}^n \xi_k(x) = n^{-1/2} S(n; x), \eqno(3.1)
$$

$$
\Gamma(L) = \Gamma_{\vec{p}}(L) = 2  \ K_R(L) \cdot \sqrt[\overline{p}]{\left[ \int_X {\bf E} |\xi(x)|^L \ \mu(dx) \right]^{\overline{p}/L}}, \eqno(3.2)
$$

$$
L_0 = L_0(\vec{p}) = \sup \{L: \ \Gamma_p(L) < \infty \}, \eqno(3.3)
$$

$$
 \gamma(z) = \gamma_{\vec{p}}(z): = \inf_{ L \in (\overline{p}, L_0)} \left\{  \frac{\Gamma^L_{\overline{p}}(L)}{z^L}  \right\}, \ z > 1. \eqno(3.4)
$$

 Let the partition $ \Delta $ be from the set  $ Y(w), \ w > 1 $ be a given. {\it Define again an important function}

 $$
 F(u) = F_{\xi}(\Delta, \vec{p},v; u) \stackrel{def}{=} \sum_{k=1}^{\infty} \gamma_{\vec{p}}(u \ v(A(k)) /w), \eqno(3.5)
 $$
if of course  $ F(u) $ there exists and $ \lim_{u \to \infty} F(u) = 0. $ \par

\newpage

{\bf Theorem 3.1. }\\

$$
Q^{(\xi)}(u) \le  F_{\xi}(\Delta, \vec{p},v; u), \ u \ge e. \eqno(3.6)
$$

{\it As a slight consequence:}

$$
Q^{(\xi)}(u) \le \inf_{w > 1} \ \inf_{\Delta \in Y(w) } \  F_{\xi}(\Delta, \vec{p},v; u), \ u \ge e. \eqno(3.7)
$$

\vspace{3mm}

We preface  the following auxiliary proposition which has by our opinion of independent interest.\\

\vspace{3mm}

{\bf Theorem 3.2.} {\it  Let $ m = \const \ge 1  $ be (not necessary to be integer) number for which }

$$
 \left[ {\bf E} |\xi(\vec{x})|^{m \ \overline{p}} \right]^{1/(m \ \overline{p}) } \in L_{ \vec{p} }. \eqno(3.8)
$$

{\it Then }

$$
\sup_n  \left[ {\bf E}  | \Sigma_n(\cdot)|_{\vec{p}}^{ m \ \overline{p}} \right]^{ 1/(m \ \overline{p}) }  \le K_R(\overline{p} \ m)  \times
  \left| \left[ {\bf E} |\xi(\vec{x})|^{m \ \overline{p}} \right]^{1/(m \overline{p})} \right|_{\vec{p}}.\eqno(3.9)
$$

\vspace{3mm}

{\bf Proof  of theorem 3.2.} \\

\vspace{3mm}

{\bf 1. Auxiliary fact: permutation inequality.} \par

\vspace{3mm}

 We will use the so-called {\it permutation inequality} in the terminology of an article \cite{Adams1};  see also \cite{Besov1}, chapter 1,  p. 24 - 26.
Indeed, let $ (Z, B, \nu) $  be another measurable space  and $ \phi:  (\vec{X},Z) = \vec{X} \otimes Z    \to R $ be measurable function. In what follows
$  \vec{X} = \otimes_k X_k.   $ Let also $  r = \const \ge \overline{p}. $ It is true the following inequality (in our notations):

$$
|\phi|_{\vec{p}, \vec{X}; r, Z } \le |\phi|_{r, Z;  \vec{p}, \vec{X}}. \eqno(3.10)
$$
 In what follows $  Z = \Omega, \ \nu = {\bf P}.  $

\vspace{3mm}

{\bf 2. Auxiliary inequality.} \par

\vspace{3mm}

 It follows from permutation inequality (3.10) that

 $$
 \sqrt[m \overline{p}] {{\bf E} |\xi|_{\vec{p}}^{m \overline{p}}}  \le
 \left| \sqrt[m \overline{p}] { {\bf E} |\xi|^{m \overline{p}}} \right|_{\vec{p}}, \ m = \const \ge 1. \eqno(3.11)
 $$

\vspace{3mm}
{\bf 3. } We deduce applying  the inequality (3.11) for the random field $ \Sigma_n(x,t)  $ and using the Rosenthal's inequality:

$$
 \sqrt[m \overline{p}] {{\bf E} |\Sigma_n|_{\vec{p}}^{m \overline{p}}}  \le K_R( m \ \overline{p} ) \cdot
 \left| \sqrt[m \overline{p}] { {\bf E} |\xi|^{m \overline{p}}} \right|_{\vec{p}}, \ m = \const \ge 1, \eqno(3.12)
 $$
which is equivalent to the assertion of theorem 3.2.\par

\vspace{3mm}

{\bf Proof of theorem 3.1} is is now completely analogous to ones in theorem 2.1 and may be omitted.\par

 \vspace{3mm}

 {\bf  Remark 3.1.} The examples on the acting of theorems 3.1 and 3.2 may be considered  at the same as in second section
 with replacement the power $ p $ on the power $  \overline{p}. $ \par

\vspace{4mm}

\section{ LIL in continuous - Lebesgue spaces.}

\vspace{3mm}

 {\bf 0. Definition of continuous-Lebesgue (Lebesgue-Riesz) space}  $ CL(p) = C(T, L_p(X)), p \ge 1.  $  \\

 \vspace{3mm}

 Let $ (X,A,\mu)  $ be again measurable space with sigma - finite separable measure $  \mu, \ T = \{t \} $  be metrizable
compact set. \par

We will say that the (measurable) function of two variables $ f = f(x,t), \ x \in X, \ t \in T $  belongs to the space
$ CL(p) = C(T, L_p(X)),$ where $ p = \const \ge 1, $ if the map $ t \to f(\cdot,t), \ t \in T  $ is continuous in  the
$ C(T) $ sense: \par

$$
\lim_{\epsilon \to 0+} \sup_{d(t,s) < \epsilon } \left[ \int_X |f(x,t) - f(x,s)|^p \ \mu(dx) \right]^{1/p} = 0. \eqno(4.1)
$$

 The norm of the function $ f(\cdot,\cdot) $ in this space is defined as follows:

$$
 ||f(\cdot, \cdot)|| C(T, L_p(X)) = ||f(\cdot, \cdot)|| CL(p) = \sup_{t \in T} |f(t, \cdot)|_p.\eqno(4.2)
$$

 These spaces  are complete separable Banach function spaces. The detail investigation of these spaces see, e.g. in
a monograph \cite{Kufner1}, p. 113 - 119.\par
  They are used, for instance, in the theory of non-linear evolution Partial Differential  Equations, see
\cite{Fujita1},  \cite{Kato1},  \cite{Kato2}, \cite{Lions1}, \cite{Lions2}, \cite{Taylor1}. \par

\vspace{3mm}

{\bf 1. Additional construction and conditions.} \par

\vspace{3mm}

{\it Assume in addition  that  $  p \ge 2 $ and that the r.f. } $  \xi(x,t),  $ {\it and with it the r.f. }
$  \xi_i(x,t)  $ {\it  are mean zero: } $ {\bf E} \xi_i(x,t) = 0, $ and denote

$$
\Sigma_n(x,t) = n^{-1/2} \sum_{i=1}^n \xi_i(x,t), \  S_n (x,t) = \sum_{i=1}^n \xi_i(x,t) = n^{1/2} \Sigma_n(x,t), \eqno(4.3)
$$

$$
\tau_p^{(n)}(t) = \int_X | \Sigma_n(x,t)|^p \ \mu(dx) = |\Sigma_n(\cdot,t)|^p_{p,X}. \eqno(4.4)
$$

\vspace{3mm}

 We intend as  before to  estimate  uniformly over numbers of summand $ n $  first of all in this section
the moments of the random variable

$$
\zeta_p^{(n)} = \zeta_p := \sup_{t \in T} \tau_p^{(n)}(t), \eqno(4.5)
$$
i.e.  the values

$$
\delta_p(Z) = \sup_n \left| \ \sup_{t \in T}  \tau_p^{(n)}(t) \  \right|_{Z, \Omega}.
$$

 Note that

 $$
 \left[ \delta_p(Z) \right]^{1/p} =  \sup_n | \Sigma_n(\cdot,\cdot)|_{p,X; \infty,T; Z,\Omega}.
 $$

\vspace{3mm}

 Some new notations:  $ \rho_{v,x}(t,s) : = $

 $$
 | \ \xi(x,t) - \xi(x,s) \ |_{v,\Omega} =  \left\{ \left[  {\bf E} |\xi(x,t) - \xi(x,s)| \right]^v    \right\}^{1/v},
\ v = \const \ge 1;
 $$

$$
W_{\gamma}(x) = \sup_{t \in T} | \ \xi(x,t) \ |_{\gamma, \Omega} = \sup_{t \in T} \left\{ {\bf E} |\xi(x,t)|^{\gamma} \right\}^{1/\gamma}, \eqno(4.6)
$$

$$
J(t,s; p,Z; \alpha,\beta) = \int_X W^{p-1}_{ (p-1) \beta Z }(x) \ \rho_{ \alpha Z,x }(t,s) \ \mu(dx), \ \alpha,\beta > 1,
1/\alpha + 1/\beta = 1;
$$

$$
r_{p,Z}(t,s) = 2 \ p \ \inf_{\alpha, \beta} \left[ K_R(\alpha Z) \ K_R^{p-1}((p-1) \beta Z) \ J(t,s; p,Z; \alpha,\beta)  \right]. \eqno(4.7)
$$

 Evidently, $ r_{p,Z}(t,s)  $ is the distance as the function on $ (t,s),  $ if it is finite.
 The minimum in the right - hand side (4.7) is calculated over all the values $ (\alpha, \beta) $ for which
 $  \alpha, \beta > 1, \ 1/\alpha + 1/\beta = 1.  $ \par

Denote also  for arbitrary  set $  T $ equipped with  distance  $ d = d(t,s) $ by
 $  N(T, d, \epsilon) $ the so-called {\it covering number: }  the minimal number of closed balls of a radii $ \epsilon > 0 $
 in the distance $  d  $ covering the set $ T. $ Obviously, $ \forall \epsilon > 0  \ N(T, d, \epsilon) < \infty $ if and
 only if the set $  T $ is precompact set relative the distance $  d. $ \par

  Further, define

$$
 \overline{\sigma}_{p,Z} \stackrel{def}{=} \sup_{t \in T} \int_X \left[ {\bf E} |\xi(x,t)|^{p Z} \right]^{1/Z}  \ \mu(dx),
\ \hat{\sigma}_{p,Z}:= K_R^p(p \ Z) \ \overline{\sigma}_{p,Z}, \eqno(4.8)
 $$

$$
\hat{r}_{p,Z}(t,s) := r_{p,Z}(t,s)/\hat{\sigma}_{p,Z},\eqno(4.9)
$$

$$
 \nu_p^p(Z) \stackrel{def}{=}   \hat{\sigma}_{p,Z} \cdot \inf_{\theta \in (0,1)}
\left[  \sum_{k=1}^{\infty} \theta^{k-1} N^{1/Z} \left(T, \hat{r}_{p,Z}, (\theta \ \hat{\sigma}_{p,Z})^k \right) \right]. \eqno(4.10)
 $$

\vspace{3mm}

{\bf 2. Theorem 4.1.} {\it If for some}   $ Z = \const \ge 1 \ \Rightarrow  \nu_p(Z) < \infty,  $  {\it then    }

\vspace{3mm}

$$
 \sup_n \left\{ {\bf E} |\Sigma_n(\cdot, \cdot)|^{pZ}_{p,\infty} \right\}^{1/p Z} \le \nu_p(Z). \eqno(4.11)
$$

\vspace{4mm}

{\bf Proof of theorem 4.1.}\\

{\bf A.} We need first of all to obtain the estimate (4.4). We have
using the Rosenthal's constants and the Minkowsky  inequality:

$$
| \ \tau_p^{(n)}(t) \ |_{Z, \Omega}  = \left| \ \int_X |\Sigma_n(x,t)|^p \ \mu(dx) \ \right|_{Z, \Omega} \le
$$

$$
 \int_X  \left| \ |\Sigma_n(x,t)|^p  \ \right|_{Z, \Omega}   \ \mu(dx)  \le
 \int_X K_R^p(p \ Z) \ \left| \ |\xi(x,t)|^p  \ \right| _{Z, \Omega}  \ \mu(dx)\le
$$

$$
K_R^p(p \ Z) \ \overline{\sigma}_{p, Z} =  \hat{\sigma}_{p, Z}. \eqno(4.12)
$$

\vspace{3mm}

{\bf B.} The estimation of a difference
$$
 \Delta \tau(t,s) =  \tau_p^{(n)}(t) - \tau_p^{(n)}(s)
$$
is more complicated. We have consequently:

$$
\Delta \tau = \int_X \left[ |\Sigma_n(x,t)|^p - |\Sigma_n(x,s)|^p  \right]  \ \mu(dx),
$$

$$
|\Delta \tau |_{Z,\Omega} \le \int_X \left| \ |\Sigma_n(x,t)|^p - |\Sigma_n(x,s)|^p \  \right|_{Z, \Omega}  \ \mu(dx) =
$$

$$
\int_X \left[ {\bf E} \left| \ |\Sigma_n(x,t)|^p -  |\Sigma_n(x,s)|^p  \ \right|^Z  \right]^{1/Z} \ \mu(dx). \eqno(4.13)
$$

  We exploit the following  elementary inequality:

$$
| \ |x|^p - |y|^p \ | \le p \cdot |x-y| \cdot \left[ |x|^{p-1} + |y|^{p-1} \right], \ x,y \in R, \eqno(4.14)
$$
and obtain after substituting into (4.12), where $  x = \Sigma_n(x,t), \ y = \Sigma_n(x,s):  \ |\Delta \tau |_{Z,\Omega}/p \le  $

$$
 \int_X \left| \ |\Sigma_n(x,t) - \Sigma_n(x,s)| \ \cdot
 \ \left[ \ |\Sigma_n(x,t)|^{p-1} + |\Sigma_n(x,s)|^{p-1} \right] \ \right|_{Z, \Omega}  \ \mu(dx). \eqno(4.15)
$$

  It follows from the H\"older's inequality

$$
|\eta_1 \eta_2|_{Z,\Omega} \le |\eta_1|_{\alpha Z, \Omega} \cdot |\eta_2|_{\beta Z, \Omega},
$$
where as before $ \alpha, \beta > 1, \ 1/\alpha + 1/\beta = 1. $ Therefore

$$
|\Delta \tau |_{Z,\Omega}/p \le \int_X \delta_1(t,s; \alpha, x) \cdot \delta_2(\beta,x) \ \mu(dx),
$$
where

$$
\delta_1(t,s; \alpha, x) =  \delta_1(t,s; \alpha, x; \Omega)  = |  \Sigma_n(x,t) - \Sigma_n(x,s) |_{\alpha Z; \Omega}, \eqno(4.16)
$$
and

$$
\delta_2(\beta,x) = \delta_2(\beta,x; p,\Omega) = \sup_{t,s \in T} \left| \ \left[|\Sigma_n(x,t)|^{p-1} + |\Sigma_n(x,s)|^{p-1} \right] \ \right|_{\beta Z, \Omega}. \eqno(4.17)
$$

 We estimate $ \delta_1(\cdot) $ using the Rosenthal's inequality:

 $$
 \delta_1(t,s; \alpha, x; \Omega) \le K_R(\alpha Z) \ | \xi(x,s) - \xi(x,s) |_{\alpha Z}= \rho_{\alpha Z,x}(t,s). \eqno(4.18)
 $$
 Further,

 $$
 \delta_2(\beta,x; p,\Omega)  \le 2 \ K_R^{p-1}(\beta (p-1))  \sup_{t \in T}  | \xi(x,t)|^{p-1}_{\beta(p-1), \Omega} =
 $$

$$
 2 \ K_R^{p-1}(\beta (p-1))\ W_{\beta(p-1)}^{p-1}(x). \eqno(4.19)
$$
 We get  after substituting into (4.16) and (4.17)

$$
|\Delta \tau |_{Z,\Omega} \le r_{p,Z}(t,s). \eqno(4.20)
$$

{\bf C.}  The proposition of theorem 4.1. follows from the  main result  of article G.Pisier \cite{Pisier1}; see also \cite{Ostrovsky504},
\cite{Ostrovsky404}, \cite{Ostrovsky405}.\par

\vspace{3mm}

{\bf Remark 4.1.}  Perhaps, it may be used in this section the so - called method of "majorizing measures", or equally "generic
 chaining"; see e.g. \cite{Fernique1}, \cite{Talagrand1}, \cite{Talagrand2}, \cite{Talagrand3}, \cite{Bednorz1}, \cite{Bednorz2},
\cite{Ostrovsky503}, \cite{Ostrovsky504}.\par
 But, by our opinion, the offered here way is more convenient for declared aims. \par

\vspace{3mm}

{\bf 3. We are ready now to formulate the main result of this section. }   Some new notations:

$$
L_0 = L_0(\vec{p}) = \sup \{L: \ \nu_p(L/p) < \infty \}, \eqno(4.21)
$$

$$
 \zeta(z) = \zeta_{\vec{p}}(z): = \inf_{ L \in (p, L_0)} \left\{  \frac{\nu^L_p(L/p)}{z^L}  \right\}, \ z > 1. \eqno(4.22)
$$

 Let the partition $ \Delta $  from the set  $ Y(w), \ w > 1 $ be a given. {\it Define again an important function}

 $$
 \Theta(u) = \Theta_{\xi}(\Delta, \vec{p},v; u) \stackrel{def}{=} \sum_{k=1}^{\infty} \zeta_{\vec{p}}(u \ v(A(k)) /w), \eqno(4.23)
 $$
if of course  $ \Theta(u) $ there exists and $ \lim_{u \to \infty} \Theta(u) = 0. $ \par

\vspace{3mm}

{\bf Theorem 4.2. }\\

$$
Q^{(\xi)}(u) \le  \Theta_{\xi}(\Delta, \vec{p},v; u), \ u \ge e. \eqno(4.24)
$$

{\it As a slight consequence:}

$$
Q^{(\xi)}(u) \le \inf_{w > 1} \ \inf_{\Delta \in Y(w) } \  \Theta_{\xi}(\Delta, \vec{p},v; u), \ u \ge e. \eqno(4.25)
$$

\vspace{3mm}

The {\bf proof}  is at the same as before.\par

\vspace{3mm}

{\bf Example 4.1.} Let $  T $ be in addition closure of bounded set in the Euclidean
 space $  R^d, \ d = 1,2,\ldots $ with ordinary distance  $  || t - s||. $
Assume that

$$
\rho_{v,x}(t,s) \le B_v(x) \ ||t - s||^{l},  \ 0 < l = \const \le 1, \eqno(4.26)
$$
where

$$
\int_X W_{2(p-1)Z}^{p-1}(x) \ B_{2 Z}(x) \ \mu(dx) \le C(p) \ Z^{b}, \ b = \const \le 1, \ Z > 2d/l; \eqno(4.27)
$$
we put $ \alpha = \beta = 1. $  Therefore
$$
J(t,s; p,Z; 2,2) \le C_1(p) \ ||t-s||^l \  Z^b.
$$

 Further, suppose

$$
 \overline{\sigma}_{p,Z} \stackrel{def}{=} \sup_{t \in T} \int_X \left[ {\bf E} |\xi(x,t)|^{p Z} \right]^{1/Z}  \ \mu(dx)
 \asymp C_2(p) \ Z^b, \ Z > 2d/l, C_2(p) \in (0,\infty). \eqno(4.28)
$$
 We deduce after simple calculations

 $$
 Q^{(\xi)} (u) \le C_3(p) \exp \left( - C_4(p) \ u^{b}  \right), \ u \ge 1. \eqno(4.29)
 $$

\vspace{3mm}

\section{ Concluding remarks.}

\vspace{3mm}

{\bf A. Lower estimates.}

\vspace{3mm}

We give a lower slightly less trivial as in (1.8) estimate for $ Q_{1/2}(u). $
 Namely, for all the non - trivial distribution $  \xi $

 $$
 Q_{1/2}(u) \ge \max \left( {\bf P}(||\xi|| > u), \ \exp \left( -C u^2 \log \log u  \right) \right).
 \   u > e^e.
 $$
 It is sufficient to note  that this estimate  is true in the one-dimensional case; see \cite{Ostrovsky403},
 \cite{Ostrovsky1}, chapter 2, section 2.6. \par

\vspace{4mm}

{\bf B. LIL in mixed Sobolev's spaces.}\\

\vspace{3mm}

 The method presented here may be used by investigation of the Law of Iterated Logarithm as well as
 Central Limit Theorem in the so-called {\it mixed Sobolev's spaces} $ W^A_{\vec{p}}, $  see e.g.
\cite{Ostrovsky401}, \cite{Ostrovsky404}. \par

 In detail, let $ (Y_k, B_k, \zeta_k), \ l = 1,2,\ldots,l  $ be again
measurable spaces with separable sigma - finite measures $ \zeta_k. $
 Let $  A  $ be closed unbounded operator acting from the space $ W_{\vec{p}, \vec{X}} $ into
 the space $ W_{\vec{p}, \vec{Y}}, $  for instance, differential operator,
 Laplace's operator or its power, may be fractional, for instance:

 $$
 A[u](x,y) = \frac{D^{(\vec{q})}u(x) - D^{(\vec{q})}u(y)}{|x-y|^{\beta}}, \ x,y \in R^d, \ x \ne y, \
 $$

 $$
 \zeta(G) = \int \int_{G} \frac{dx dy}{|x-y|^{\alpha}}, \ \alpha,\beta = \const \in [0,1],  \alpha + \beta p < d, \ G \subset R^{2d},
 $$

$$
\vec{q} = \{q_1, q_2,  \ldots, q_d \}, \ q_s = 0,1, \ldots,
$$

$$
D^{(\vec{q})}u(x) = \frac{\partial^{q_1}}{\partial x_1^{q_1}} \frac{\partial^{q_2}}{\partial x_2^{q_2}} \ldots    \frac{\partial^{q_d}}{\partial x_d^{q_d}}u(x).
$$
  The norm in this space may be defined as follows (up to closure):

 $$
|f|W^A_{\vec{p}}\stackrel{def}{=} \max \left( |f|_{\vec{p}}, |Af|_{\vec{p}} \right).  \eqno(5.1)
 $$

 Analogously to the proof of theorems 4.1, 3.2 may be obtained the following results.\par

\vspace{3mm}

{\bf Theorem 5.1.} {\it  Let $ m = \const \ge 1  $ be (not necessary to be integer) number for which }

$$
 \left[ {\bf E} |\xi(\vec{x})|^{m \ \overline{p}} \right]^{1/(m \ \overline{p}) } \in L_{ \vec{p} }, \hspace{5mm}
 \left[ {\bf E} |A[\xi](\vec{x})|^{m \ \overline{p}} \right]^{1/(m \ \overline{p}) } \in L_{ \vec{p} }. \eqno(5.2)
$$

{\it Then }

$$
\sup_n \left[ {\bf E} |A[ \Sigma_n](\cdot)|_{\vec{p}}^{ m \ \overline{p}} \right]^{ 1/(m \ \overline{p}) }  \le K_R(\overline{p} \ m)  \times
  \left| \left[ {\bf E} |A[\xi](\vec{x})|^{m \ \overline{p}} \right]^{1/(m \overline{p})} \right|_{\vec{p}},\eqno(5.3)
$$
{\it  Hence  the propositions of theorem 2.1, 3.1 remains true  in this norm. } \par

\vspace{3mm}

{\bf B. LIL for dependent r.v. in anisotropic spaces.}\\

\vspace{3mm}

{\it We refuse in this section on the assumption about independence of random vectors} $ \{ \xi_k(\cdot)  \}.  $ \\

\vspace{3mm}

{\bf  Martingale case.} \\

\vspace{3mm}

 We suppose as before that  $ \{ \xi_k(\cdot)  \}  $ are mean zero and form a strictly stationary sequence with values
 in  Bochner's (mixed) space $  L_{\vec{p}}, \  \overline{p} \ge 2. $\par

  Assume in addition that $ \{ \xi_k(\cdot)  \}  $ form a martingale difference sequence
 relative certain filtration $ \{  F(k) \},  \ F(0) = \{ \emptyset, \Omega\},  $
 $$
  {\bf E} \xi_k/F(k) = \xi_k, \hspace{5mm}  {\bf E} \xi_k/F(k-1) = 0, \ k= 1,2,\ldots.
 $$

 Then the propositions  of theorems 2.1, 3.1 remains true; the estimate  of theorem 3.2 is also true  up to multiplicative
absolute constant. \par

 Actually, the Law of Iterated Logarithm  for one-dimensional  martingales
 see in the classical monograph of  Hall P., Heyde C.C.
  \cite{Hall1}, chapter 2; the Rosenthal's constant for the sums  of martingale differences with at the same up to
  multiplicative constant coefficient is obtained by A.Osekowski \cite{Osekowski1}, \cite{Osekowski2}.
 See also \cite{Ostrovsky3}.\par

\vspace{3mm}

{\bf  Mixingale case.} \\

\vspace{3mm}

 We suppose again  that  $ \{ \xi_k(\cdot)  \}  $ are mean zero and form a strictly stationary sequence,
 $  \overline{p} \ge 2. $  This  sequence is said to be {\it mixingale, } in the terminology of the book
\cite{Hall1}, if it satisfies this or that mixing condition.\par

 We  consider here only the  superstrong mixingale. Recall that the superstrong, or $  \beta = \beta(F_1, F_2) $
index between two sigma - algebras  is defined as follows:

$$
\beta(F_1, F_2) = \sup_{A \in F_1, B \in F_2, {\bf P}(A) {\bf P}(B) > 0 } \left| \frac{{\bf P}(AB) - {\bf P}(A) {\bf P}(B)}{{\bf P}(A) {\bf P}(B)} \right|.
$$

 Denote

 $$
 F_{-\infty}^0  = \sigma(\xi_s, \ s \le 0),  \hspace{5mm} F_n^{\infty} = \sigma(\xi_s, \ s \ge n), \eqno(5.4)
 $$

$$
\beta(n) = \beta \left(F_{-\infty}^0 , F_n^{\infty} \right),
$$

 The  sequence $ \{\xi_k \} $ is said to be {\it superstrong mixingale, } if $ \lim_{n \to \infty} \beta(n) = 0. $ \par
This notion  with some applications was introduced and investigated by B.S.Nachapetyan and R.Filips \cite{Nachapetyan1}.
See also   \cite{Ostrovsky3}, \cite{Ostrovsky1}, p. 84 - 90. \par

\vspace{3mm}

 Introduce the so-called mixingale Rosenthal coefficient:

 $$
 K_M(m) =  m \ \left[  \sum_{k=1}^{\infty} \beta(k) \ (k+1)^{ (m -2)/2  }   \right]^{1/m}, \ m \ge 1.\eqno(5.5)
 $$

 B.S.Nachapetyan in \cite{Nachapetyan1} proved that for the superstrong centered  strong stationary strong mixingale
sequence $  \{ \eta_k \} $ with $ K_M(m) < \infty $ the following estimate is true:

$$
\sup_{n \ge 1}  \left| n^{-1/2} \sum_{k=1}^n \eta_k \right|_m \le C \cdot K_M(m) \cdot |\eta_1|_m, \eqno(5.6)
$$
so that the "constant" $ K_M(m) $ play at the same role for mixingale as the Rosenthal  constant  for independent variables.\par

 As a consequence:  theorems 3.1 (and 3.2)  remains true for the strong mixingale sequence $ \{  \xi_k \}: $
theorem 3.1 under conditions: $ K_M(m \overline{p}) < \infty, $
with replacing $  K_R( m \ \overline{p} ) $ on the expression $  K_M( m \ \overline{p} ). $  \par

\end{document}